\documentclass[journal]{IEEEtran}
\usepackage{amsmath,amssymb,graphicx,paralist,subfigure,pdfpages,cite,algorithm,algpseudocode,paralist,amsthm}
\usepackage{verbatim}
\usepackage{setspace}
\newtheorem{lem}{Lemma} 
\newtheorem{theorem}{Theorem}

\newtheorem{rmk}{Remark}
\newtheorem{assump}{Assumption}

\def\bpi{\boldsymbol\pi}

\def\mc{\mathcal}
\def\mb{\mathbf}
\def\mbb{\mathbb}
\def\ra{\rightarrow}
\def\la{\leftarrow}
\def\P{\mathbf{P}}

\IEEEoverridecommandlockouts
\def\AB{\textbf{\texttt{AB}}}
\def\SAB{\textbf{\texttt{S-AB}}}

\def\SGP{\textbf{\texttt{SGP}}}
\def\SA{\textbf{\texttt{S-ADDOPT}}}

\def\DSGD{\textbf{\texttt{DSGD}}}
\def\DSGT{\textbf{\texttt{GT-DSGD}}}
\def\SAGA{\textbf{\texttt{SAGA}}}

\def\PS{\textbf{\texttt{Push-SAGA}}}
\def\AS{\textbf{\texttt{AB-SAGA}}}

\def\mt{\times}
\def\mbb{\mathbb}
\def\mb{\mathbf}
\def\mc{\mathcal}
\def\wh{\widehat}

\def\ol{\overline}
\def\ul{\underline}

\def\bds{\boldsymbol}
\newcommand{\mn}[1]{{\left\vert\kern-0.25ex\left\vert\kern-0.25ex\left\vert\kern0.3ex #1 
\kern0.3ex\right\vert\kern-0.25ex\right\vert\kern-0.25ex\right\vert}}

\begin{document}
\title{\huge \textbf{Variance reduced stochastic optimization over directed graphs with row and column stochastic weights}}
\author{
Muhammad I. Qureshi, Ran Xin, Soummya Kar, and Usman A. Khan
\thanks{MIQ and UAK are with the Electrical and Computer Engineering Department (ECE) at Tufts university, Medofrd, MA, USA; \texttt{muhammad.qureshi@tufts.edu}, \texttt{khan@ece.tufts.edu}. RX and SK are with the ECE Department at Carnegie Mellon University, Pittsburgh, PA, USA; \texttt{\{ranx,soummyak\}@andrew.cmu.edu}. The authors acknowledge the support of NSF  under awards  CMMI-1903972 and CBET-1935555.}
}
\maketitle

\begin{abstract}
This paper proposes~\AS, a first-order distributed stochastic optimization method to minimize a finite-sum of smooth and strongly convex 
functions distributed over an arbitrary directed graph. $\AS$ removes the uncertainty caused by the stochastic gradients using a node-level variance reduction and subsequently employs network-level gradient tracking to address the data dissimilarity across the nodes. Unlike existing methods that use the nonlinear push-sum correction to cancel the imbalance caused by the directed communication, the consensus updates in~$\AS$ are linear and uses both row and column stochastic weights. We show that for a constant step-size,~$\AS$ converges linearly to the global optimal. We quantify the directed nature of the underlying graph using an explicit directivity constant and characterize the regimes in which~$\AS$ achieves a linear speed-up over its centralized counterpart. Numerical experiments illustrate the convergence of $\AS$ for strongly convex and nonconvex problems.


\end{abstract}

\section{Introduction} \label{intro}
Stochastic optimization is relevant in many signal processing, machine learning, and control applications~\cite{bottou2018optimization,BigData_Bajwa,MRD_camsap_2011,nc1,li2022communication}. In very large-scale problems, data is usually geographically distributed making centralized methods practically infeasible. Distributed solutions are thus preferable where individual nodes perform local updates with the help of data fusion among the nearby nodes~\cite{GT_SAGA_SPM,AB_proc,math10030357,wang2021decentralized,SONG2022110096}. The problem of interest can be written as 
\[
\P:\quad \min_{\mb x\in\mbb R^p} F(\mb{x}) := \frac{1}{n} \sum_{i=1}^{n} f_i(\mb{x}), ~~ f_i(\mb{x}) := \frac{1}{m_i} \sum_{j=1}^{m_i} f_{i,j}(\mb{x}),
\]
where each local cost function~${f_i: \mbb{R}^p \ra \mbb{R}}$ is private to node~$i$ and is further decomposed into~$m_i$ component cost functions~${\{f_{i,j}:\mbb{R}^p \ra \mbb{R}\}_{j=1}^{m_i}}$. When the underlying optimization problem is smooth and strongly convex, the goal is to find the unique minimizer~$\mb x^*$ of the global cost~$F$, assuming that the network consists of~$n$ nodes communicating over an arbitrary strongly-connected directed graph.

Distributed first-order stochastic methods for problem~$\P$ are well studied in the literature. Early work includes~\cite{DSGD_nedich,distributed_Jianshu} that is applicable to undirected graphs. Stochastic gradient push ($\SGP$~\cite{SGP_nedich,SGP_ICML,SGP_AsyncOlshevsky}) extends $\DSGD$ (distributed stochastic gradient descent~\cite{DSGD_nedich}) to directed graphs using push-sum consensus~\cite{ac_directed0,pushsum_Hadjicostis}. Both $\DSGD$ and $\SGP$ suffer from a steady-state error caused by the difference in global and local cost functions, i.e.,~${\|\nabla F(\mb{x}^*) - \nabla f_i(\mb{x}^*)\|}$, and the variance introduced by the stochastic gradients. Over arbitrary directed graphs, $\SA$ \cite{saddopt} compensates for the heterogeneity of local cost functions with the help of gradient tracking~\cite{di2016next,AugDGM,harnessing,add-opt}. However, the steady-state error remains in effect due to the variance. A recent work $\PS$~\cite{pushsaga} benefits from a variance reduction technique~\cite{SAGA} to eliminate the uncertainty caused by the stochastic gradients. Both $\SA$ and $\PS$ use push-sum correction to implement consensus nonlinearly and divide by the estimates of the right Perron eigenvector of the underlying column stochastic weight matrix. Such correction is not required when the weights are doubly stochastic as is the case over undirected (or weight-balanced directed) graphs; see~\cite{EXTRA,harnessing,Undirected_Olshev,DSGT_Pu,GT_SAGA,GT_SAGA_ncvx,GT_SARAH_ncvx} for related work. 

In this paper, we present~$\AS$, a first-order distributed stochastic optimization method that is applicable to arbitrary directed graphs; see also~\cite{AB,ABN}. Similar to the methods in~\cite{GT_SAGA,pushsaga}, $\AS$ eliminates the uncertainty caused by the stochastic gradients with the help of variance reduction and addresses the global vs. local cost gaps, due to data dissimilarity across different nodes, using gradient tracking. Unlike \PS~\cite{pushsaga}, however, $\AS$ uses both row and column stochastic weights to ensure consensus, thus eliminating the need of estimating the Perron eigenvector required in push-sum methods; see~\cite{AB} for the~$\AB$ algorithm. The main contributions of this paper are summarized next: \begin{inparaenum}[(i)]
\item We demonstrate the linear convergence of~$\AS$ to the global optimizer~$\mb{x}^*$ of smooth and strongly convex problems;
\item We quantify the performance of~$\AS$ over directed graphs and encapsulate the directed nature of the communication in a directivity constant~${\psi \geq 1}$, which is unity for undirected graphs;
\item We provide explicit expressions for the gradient computation and communication complexities, and show that $\AS$ achieves linear speedup over its centralized counterpart \SAGA~\cite{SAGA}.
\end{inparaenum}

We now describe the rest of the paper. Section~\ref{algo_dev} provides the algorithm development and formally describes~$\AS$. Section~\ref{main_res} describes the assumptions and the main convergence results, whereas Section~\ref{conv} provides the detailed convergence analysis. Finally, Section~\ref{sims} presents the numerical experiments on strongly convex and nonconvex problems, and Section~\ref{conc} concludes the paper.

\textbf{Basic Notation:} We use upper case letters to represent matrices and lower case bold letters for vectors. We define~$I_n$ as~${n \mt n}$ identity matrix and~$\mb{1}_n$ as a column vector of~$n$ ones. From Perron Frobenius theorem~\cite{hornjohnson}, for a primitive row stochastic matrix~${\ul A \in \mbb{R}^{n \mt n}}$ (column stochastic matric~${\ul B \in \mbb{R}^{n \mt n}}$), we define~${\ul{A}^\infty := \lim_{k \ra \infty}\ul{A}^k = \bpi_r^\top \mb{1}_n}$ (and ${\ul{B}^\infty := \lim_{k \ra \infty}\ul{B}^k = \mb{1}_n^\top \bpi_c}$), where~$\bpi_r$ is the left eigenvector of~$\ul{A}$ ($\bpi_c$ is the right eigenvector of~$\ul{B}$), corresponding to the unique eigenvalue~$1$. We further denote the largest element of a vector~$\bpi_r$ as~$\ol{\pi}_r$ and the smallest element as~$\ul{\pi}_r$, and define the ratios~${h_r := \ol{\pi}_r/ \ul{\pi}_r}$ and~${h_c := \ol{\pi}_c/ \ul{\pi}_c}$. We next define the spectral radius of matrix~$\ul{A}$ as~$\rho(\ul{A})$. We denote~${\|\cdot\|_2}$ as the Euclidean norm and~${\mn{\cdot}}$ as the matrix norm. Since~${\rho(\ul{A} - \ul{A}^\infty)<1}$ and~${\rho(\ul{B} - \ul{B}^\infty)<1}$, it can be shown that there exist matrix norms~$\mn{\cdot}_{\bpi_r}$ and~$\mn{\cdot}_{\bpi_c}$, formally defined in~\cite{khan_cdc1:19}, such that~${\sigma_A:=\mn{\ul{A} - \ul{A}^\infty}_{\bpi_r}}$ and~${\sigma_B:=\mn{\ul{B} - \ul{B}^\infty}_{\bpi_c}}$. 

\section{Algorithm development} \label{algo_dev}
We motivate the proposed algorithm with the help of a recent work $\DSGT$~\cite{DSGT_Pu, DSGT_Xin_ncvx}, which adds gradient tracking to the well known $\DSGD$ \cite{DSGD_nedich}. The $\DSGT$ algorithm can be described as follows. Let ${\ul{W} = \{\ul{w}_{ij} \}}$ be the network weight matrix such that~${\ul{w}_{ij} \neq 0}$, if and only if node~$i$ can receive information from node~$j$. Let~${\mb{x}_i^k, \mb{w}_i^k}$, both in~$\mbb{R}^p$ be the state vectors at each node~$i$ and iteration~$k$. Then ${\forall k \geq 0}$, $\DSGT$ at each node~$i$ is given by
\begin{align*}
    \mb{x}^{k+1}_i &= \sum_{r=1}^n \ul{w}_{i r} \mb{x}^k_i - \alpha \mb{w}^{k}_i,\\
    \mb{w}^{k+1}_i &= \sum_{r=1}^n \ul{w}_{i r} \mb{w}^k_i + \nabla f_{i,s_i^{k+1}} \left(\mb{x}^{k+1}_i \right) - \nabla f_{i,s_i^{k}} \left(\mb{x}^k_i \right),
\end{align*}
where~$s_i^k$ is an index drawn uniformly at random from the index set $\{1,\cdots,m_i\}$ and~$\nabla f_{i,s_i^{k}}(\mb{x}^k_i)$ is the gradient of the~$s_i^{k}$-th component cost function~$f_{i,s_i^k}$ (and not the full local gradient~$\nabla f_i$). The~$\mb w^k_i$-update in $\DSGT$ is based on dynamic average consensus~\cite{DAC} and essentially tracks the global gradient~$\nabla F$, asymptotically, see~\cite{di2016next,AugDGM,harnessing,add-opt} for more details. The~$\mb x^k_i$-update consequently implements a descent in the global gradient direction~$\mb w^k_i$. Assuming that the variance of local stochastic gradients is bounded, i.e.,~${\mathbb{E}_{s^k_i}[\|\nabla f_{i,s^k_i}(\mb x^k_i)-\nabla f_i(\mb x_i^k)\|_2^2 ~|\:\mb x_i^k] \leq\sigma^2}$, and the global cost is~$\ell$-smooth and $\mu$-strongly convex, $\DSGT$ converges linearly to the neighborhood of the optimal solution, i.e.,
\begin{equation*}\label{DSGT_convergence}
\limsup_{k\rightarrow\infty}\frac{1}{n}\sum_{i=1}^{n} \mbb{E}[\|\mb x^k_i-\mb x^*\|_2^2]
=  \mc{O}\left(\frac{\alpha}{n\mu}\:\sigma^2
+ \frac{\alpha^2\kappa^2}{(1-\lambda)^3}\:\sigma^2\right),
\end{equation*}
for a sufficiently small constant stepsize~$\alpha$, where~$(1-\lambda)$ is the spectral gap of~$\ul{W}$ and~$\kappa$ is the condition number of $F$. We note that the steady state error in $\DSGT$ depends on the variance of the stochastic gradients~$\sigma^2$. Moreover,~$\DSGT$ is applicable to undirected graphs since it requires the weight matrix~$\ul{W}$ to be doubly stochastic. 

In this paper, we propose $\AS$ that removes the steady state error in $\DSGT$ with the help of a variance reduction technique based on the~$\SAGA$ method~\cite{SAGA}. Moreover, $\AS$ is applicable to arbitrary directed graphs as it only requires a row stochastic matrix~$\ul A$ and column stochastic matrix~$\ul B$. The complete implementation details are formally described in Algorithm~\ref{algo}. We note that for each~$\mb{x}_i^k$ update,~$\AS$ requires~${c \in \mbb{N}}$ communication rounds and for each~$\mb{w}_i^k$ update, it requires~${d \in \mbb{N}}$ communication rounds. For ease of notation, we write the~$ij$-th element of~$\ul{A}^c$ as~$\{a_{ij}\}$ and~$\ul{B}^d$ as~$\{b_{ij}\}$, for some~${c,d \in \mbb N}$ formally defined later. Each node~$i$ updates~${\mb{x}_i^k}$ which estimates the global minimum~$\mb{x}^*$, and~${\mb{w}_i^k}$ which tracks the gradient~$\nabla F(\mb{x}_i^k)$ of the global cost~$F$ using~$\SAGA$-based local gradient update~$\mb{g}_i^{k+1}$. We remark that each node~$i$ requires additional storage~$\mc{O}(p m_i)$ to maintain the gradient table~$\{\nabla f_{i,j}(\mb v^k_{i,j})\}_{j=1}^{m_i}$ as is standard in \SAGA-based methods. This storage cost can be reduced to~$\mc{O}(m_i)$ for certain problems~\cite{SAGA}. 
\begin{algorithm}[t] 
	\caption{~\AS~at each node~$i$} \label{algo}
	\begin{spacing}{1.1}
    \begin{algorithmic}[1] 
		\Require ${\mb{x}_i^0\in\mbb R^p},~{\mb{w}_i^0=\mb{g}_i^0=\nabla f_{i}(\mb{x}_{i}^{0})},~{\mb{v}_{i,j}^1 = \mb{x}_i^0},$
		
		~~${\forall j \in \{1,\cdots,m_i\}}, {\alpha>0},~{\{a_{ir}\}_{r=1}^n},~{\{b_{ir}\}_{r=1}^n}$,
		
		~~Gradient table:~$\{\nabla f_{i,j}(\mb{v}_{i,j}^0)\}_{j=1}^{m_i}$
		
		\For{$k = 0,1,2,\dots$}
		\State$\mb{x}_{i}^{k+1} \la \sum_{r=1}^n a_{ir} (\mb{x}_{r}^{k} - \alpha\cdot\mb{w}_{i}^{k})$ 
		\State \textbf{Select}~$s_{i}^{k+1}$ uniformly at random from~$\{1,\dots,m_i\}$
		\State $\mb{g}_{i}^{k+1} \la \nabla f_{i,s_i^{k+1}}(\mb{x}_{i}^{k+1}) - \nabla f_{i,s_i^{k+1}}(\mb{v}_{i,s_i^{k+1}}^{k+1}) + \frac{1}{m_i}\sum_{j=1}^{m_i}\nabla f_{i,j}(\mb{v}_{i,j}^{k+1})$ 
		\State \textbf{Replace}~$\nabla f_{i,s_i^{k+1}}(\mb{v}_{i,s_i^{k+1}}^{k+1})$ by~$\nabla f_{i,s_i^{k+1}}(\mb{x}_{i}^{k+1})$ in the gradient table
		\State $\mb{w}_{i}^{k+1} \la \sum_{r =1 }^{n}b_{ir} (\mb{w}_{r}^{k} + \mb{g}_i^{k+1} - \mb{g}_i^{k})$
		\If{$j = s_{i}^{k+1}$,}
		$\mb{v}_{i,j}^{k+2} \la \mb{x}_i^{k+1}$,~\textbf{else}~$\mb{v}_{i,j}^{k+2} \la \mb{v}_{i,j}^{k+1}$
		\EndIf
		\EndFor
	\end{algorithmic}
	\end{spacing}
\end{algorithm}

\vspace{-0.2cm}
\section{Assumptions and Main Results} \label{main_res}
We first describe the assumptions below.
\begin{assump}\label{weights} The network of nodes communicate over a strongly connected arbitrary directed graph.
\end{assump}
\vspace{-0.3cm}
\begin{assump}\label{smooth_convex} The global cost function $F$ is~${\mu}$-strongly convex and each component cost~$f_{i,j}$ is~${\ell}$-smooth.
\end{assump}
Assumption 1 ensures that the resulting weight matrices ${\ul A=\{\ul a_{ir}\}}$ and~${\ul B=\{\ul b_{ir}\}}$ are both irreducible and primitive. These requirements can be fulfilled if each node~$i$ has the knowledge of its in-degree~$d_i^{\text{in}}$ and its out-degree~$d_i^{\text{out}}$. Then the weights can be locally chosen as~${\ul{a}_{ir} = 1/d_i^{\text{in}}}$~for each incoming neighbour~$r$,~and~${\ul{b}_{ir} = 1/d_i^{\text{out}}}$~for each outgoing neighbour~$r$. Next, Assumption 2 ensures that the global cost~$F$ is~$\ell$-smooth and~$\mu$-strongly convex and therefore has a unique minimizer~$\mb{x}^*$. We note that the local cost functions~$f_i$'s are not necessarily strongly convex, which is a relaxed condition than the one for~$\PS$. Based on these assumptions, we now present the main results.

\begin{theorem} \label{th1}
Consider problem~$\P$ and let Assumptions~\ref{weights} and~\ref{smooth_convex} hold. For the step-size~${\alpha \in (0, \ol{\alpha})}$,~\AS~linearly converges to the global minimizer~$\mb{x}^*$. In particular, when ${\alpha=\ol{\alpha}}$,~$\AS$ achieves an~$\epsilon$-optimal solution in
\begin{equation*}
    \Gamma = \mc{O} \left(\max \left\{\kappa \psi, \tfrac{\kappa^2M}{m}, M\right\} \log \tfrac{1}{\epsilon}\right)
\end{equation*}
gradient computations, with~$\left( c+d \right)$~communication rounds per iteration, for all ${c=\lceil\:\ol{c}\:\rceil}$ and~${d=\lceil\:\ol{d}}\:\rceil$ such that 
\begin{align*}
    \ol{c} &:= \tfrac{\log \left( \frac{90512 n M \kappa}{m (1-\sigma_B^{2d})} \sqrt{\frac{h_r h_c}{\bpi_r^\top \bpi_c}} \right)}{\log \frac{1}{\sigma_A}}, \quad \ol{d} := \tfrac{\log \left(\frac{1265 \kappa }{\bpi_r^\top \bpi_c}\sqrt{\frac{n M h_c}{m}}\right)}{\log \frac{1}{\sigma_B}},
\end{align*}
where~${M := \max_i m_i}, {m := \min_i m_i},$ ${\kappa := \ell/\mu}$~is the condition number, and~${\psi \geq 1}$ is the directivity constant.
\end{theorem}

The formal proof of the Theorem~\ref{th1} is provided in Section~\ref{conv}. The following remarks summarize its key attributes. 
\begin{rmk}
We note that for well-connected networks, i.e., when~$\sigma_A$ and~${\sigma_B}$ are small, we have that~${\ol{c}\approx 0}$ and~${\ol{d} \approx 0}$. Thus, we get~${c=1}$ and~${d=1}$, and $\AS$ converges with a single round of communication per iteration. Furthermore, in contrast to~\cite{pushsaga}, the gradient computation complexity~$\Gamma$ is independent of the spectral gap~$(1-\sigma_A)$ and~$(1-\sigma_B)$. 
\end{rmk}

\begin{rmk}
Theorem~\ref{th1} quantifies the directed nature of the underlying graph in terms of an explicit directivity constant~${\psi:=\frac{\sqrt{h_r h_c}}{n (\bpi_r^\top \bpi_c)}}$. Clearly,~${\psi=1}$ for undirected networks; thus $\AS$ and its convergence proof are naturally applicable to undirected graphs. 
\end{rmk}

\begin{rmk}
When each node possess a large dataset such that~${M \approx m \gg \kappa^2 \psi}$, $\AS$ achieves an~$\epsilon$-optimal solution in~$\mc{O}\left(M \log \epsilon^{-1}\right)$ gradient computations per node. We note that this complexity is~$n$ times better than the centralized complexity~$\mc{O}\left(nM \log \epsilon^{-1}\right)$ of $\SAGA$~\cite{SAGA} that processes all data at a single location. \end{rmk}
\vspace{-0.25cm}
\section{Convergence of~\AS} \label{conv}
In this section, we formalize the convergence analysis. It can be verified that~$\AS$ described in Algorithm~\ref{algo} can be compactly written in vector-matrix format as
\begin{subequations}\label{ABVRv}
\begin{align}
    \mb{x}^{k+1} &= A^c (\mb{x}^{k} - \alpha \mb{w}^{k}),\\
    \mb{w}^{k+1} &= B^d (\mb{w}^{k} + \mb{g}^{k+1} - \mb{g}^{k});
\end{align}
\end{subequations}
where~$\mb{x}^{k}, \mb{w}^{k}$~and~$\mb{g}^{k}$~are the global state vectors in~$\mbb{R}^{pn}$ concatenating the local state vectors~$\mb{x}_i^{k}, \mb{w}_i^{k}$~and~$\mb{g}_i^{k}$ in~$\mbb{R}^{p}$ respectively. Similarly,~${A:= \ul{A} \otimes I_p}$ and~${B:= \ul{B} \otimes I_p}$, in~$\mbb{R}^{pn \mt pn}$, are the global weight matrices, whereas~$c$ and~$d$ denotes the communication rounds per iterate. We next define four error terms to aid the convergence analysis of~$\AS$
\begin{enumerate}[(i)]
    \item Network agreement error:~$\mbb{E} \| \mb{x}^k - A^\infty \mb{x}^k \|^2 $;
    \item Optimality gap:~${\mbb{E} \| \wh{\mb{x}}^k - \mb{x}^* \|^2}$;
    \item Mean auxiliary gap:~${\mbb{E} [ \mb{t}^k]}$;
    \item Gradient tracking error:~${\mbb{E} \| \mb{w}^k - B^\infty \mb{w}^k \|^2}$;
\end{enumerate} 
where~${\wh{\mb x}^k:= \bpi_r^\top \mb x^k}$ and~${\mb{t}^{k} := \sum_{i=1}^{n} (\frac{1}{m_i} \sum_{j=1}^{m_i} \| \mb{v}^{k}_{i,j} - \mb{x}^* \|_2^2 )}$.

In order to establish linear convergence of~\AS, we would like to show that all of the error terms described above, linearly decay to zero, eventually implying that~${\mb{x}_i^k \ra \mb{x}^*}$ for each node~$i$. We next describe the LTI system that governs the convergence rate of~$\AS$ in terms of the above error quantities in Lemma~\ref{sys_lem}.

\begin{lem} \label{sys_lem}
Consider~$\AS$ under Assumptions~\ref{weights} and~\ref{smooth_convex}. If~${\alpha \leq \min \left\{\frac{1}{35 \ell \sqrt{h_r h_c}}, \frac{\mu}{288 n \ell^2 (\bpi_r^\top \bpi_c)} \right\}}$,~${c \geq \frac{\log(4 n)}{\log(1/\sigma_A)}}$, and ${d \geq \frac{\log(4 n)}{\log(1/\sigma_B)}}$; then~${\forall k > 0, \mb{t}^{k+1} \leq G_\alpha \mb{t}^k}$, where~${\mb{t}^k \in \mbb{R}^4}$ and~${G_\alpha  \in \mbb{R}^{4 \mt 4}}$ are defined as
\begin{align}
\mb{t}^k &:=
\left[ {\begin{array}{c}
\mbb{E}[\|\mb{x}^k - A^\infty \mb{x}^k \|^2 _{\pi_r}] \\
\mbb{E}[n\|\wh{\mb{x}}^k - \mb{x}^* \|_2 ^2] \\
\mbb{E}[\mb{t}^k] \\
\mbb{E}[\ell^{-2}\|\mb{w}^k - B^\infty \mb{w}^k \|^2 _{\pi_c}]
\end{array} } \right],
\end{align}
\begin{align}
G_\alpha &:= \left[ {\begin{array}{c c c c}
\frac{3}{4}  & \alpha^2 g_1 \sigma_A^{2c} & \alpha^2 g_2 \sigma_A^{2c} & \alpha^2 g_3 \sigma_A^{2c} \\
\alpha g_{4} & 1- \alpha g_5 & \alpha^2 g_6 & \alpha g_7 \\
\frac{2}{m \ul{\pi}_r} &\frac{2}{m} & 1-\frac{1}{M} &0 \\
\frac{146 n \sigma_B^{2d}}{(1-\sigma_B^{2d}) \ul{\pi}_r \ul{\pi}_c} & \frac{97 n \sigma_B^{2d}}{(1-\sigma_B^{2d})\ul{\pi}_c} & \frac{26 \sigma_B^{2d}}{(1-\sigma_B^{2d})\ul{\pi}_c} & \frac{3}{4}
\end{array} } \right];\nonumber
\end{align}
and the constants are
\begin{align*}
\begin{array}{ll}
g_{1} := \frac{40 \ell^2 n \|\bpi_c\|^2_2 \ol{\pi}_r}{1-\sigma_A^{2c}}, \quad&g_{2} := \frac{16 \ell^2 \|\bpi_c\|^2_2 \ol{\pi}_r}{1-\sigma_A^{2c}}, \\
g_{3} := \frac{8 \ell^2 \ol{\pi}_r \ol{\pi}_c}{1-\sigma_A^{2c}},\quad&g_4 := \frac{8 \ell^2 n \bpi_r^\top \bpi_c}{\mu \ul{\pi}_r},\\
g_{5} := \frac{\mu n \bpi_r^\top \bpi_c}{4},\quad&g_6 := 3 \ell^2 n (\bpi_r^\top \bpi_c)^2,\\
g_{7} := \frac{5 \ell^2 \|\bpi_r \|^2_2 \ol{\pi}_c}{\mu \bpi_r^\top \bpi_c}.
\end{array}
\end{align*}
\end{lem}
The proof of Lemma~\ref{sys_lem} is standard and follows similar procedures as in~\cite{khan_cdc1:19,pushsaga}. With the help of this lemma, we next prove Theorem~\ref{th1} based on the following key result. 
\begin{lem} \label{spec_rad_lem}
\cite{hornjohnson} Let~${A\in \mbb{R}^{n \mt n}}$ be a non-negative matrix and ${\mb{x} \in \mbb{R}^n}$ be a positive vector. If~${A \mb{x} \leq \beta \mb{x}}$ for ${\beta > 0}$, then ${\rho(A) \leq \mn{A}_{\infty}^{\mb{x}} \leq \beta}$, where~$\mn A_\infty^{\mb x}$ is the matrix norm induced by the weighted max-norm~$\|\cdot\|_\infty^\mb{x}$ where~${\mb{x}>\mb{0}_n}$.
\end{lem}

\textbf{Proof of Theorem~\ref{th1}}: We note that the system matrix~$G_\alpha$ described in Lemma~\ref{sys_lem} is non-negative. From Lemma~\ref{spec_rad_lem}, if there exists a positive vector~${\bds{\delta}\in \mbb{R}^4}$ and a constant~$\gamma$, such that ${G_\alpha \bds{\delta} \leq \gamma \bds{\delta}}$ element-wise, then we ensure that ${\rho(G_\alpha) \leq \mn{G_\alpha}_{\infty}^{\bds{\delta}} \leq \gamma}$. To this aim, let~${\bds{\delta} := [\delta_1~\delta_2~\delta_3~\delta_4]^\top}$ has all positive elements and set~${\gamma = \left(1 - \tfrac{1}{2}\alpha g_5\right)}$ then, for~$G_\alpha$ in Lemma~\ref{sys_lem}, the following set of inequalities must hold:
\begin{align}
    &\frac{\alpha g_5}{2} + \frac{\sigma_A^{2c}}{\delta_1} \left(\alpha^2 g_1 \delta_2 + \alpha^2 g_2 \delta_3 + \alpha^2 g_3 \delta_4\right) \leq \frac{1}{4}, \label{frt}\\
    &\alpha g_6 \leq \frac{g_5}{2} \frac{\delta_2}{\delta_3} - g_4 \frac{\delta_1}{\delta_3} - g_7 \frac{\delta_4}{\delta_3}, \label{sec}\\
    &\frac{\alpha g_5}{2} \leq \frac{1}{M} - \frac{2 \ul{\pi}_r^{-1}}{m} \frac{\delta_1}{\delta_3} - \frac{2}{m} \frac{\delta_2}{\delta_3},\label{trd}\\
    &\frac{\alpha g_5}{2} \leq \frac{1}{4} - \frac{\sigma_B^{2d}}{\delta_4} \left(\frac{146 n \ul{\pi}_r^{-1} \ul{\pi}_c^{-1}}{1-\sigma_B^{2d}} \delta_1 + \frac{97 n \ul{\pi}_c^{-1}}{1-\sigma_B^{2d}} \delta_2\right) \nonumber\\
    &~~~~~~~~- \frac{\sigma_B^{2d}}{\delta_4} \left( \frac{26 \ul{\pi}_c^{-1}}{1-\sigma_B^{2d}} \delta_3 \right). \label{frth}
\end{align}

\begin{figure*}[!ht]
\centering
\includegraphics[height=1.86in]{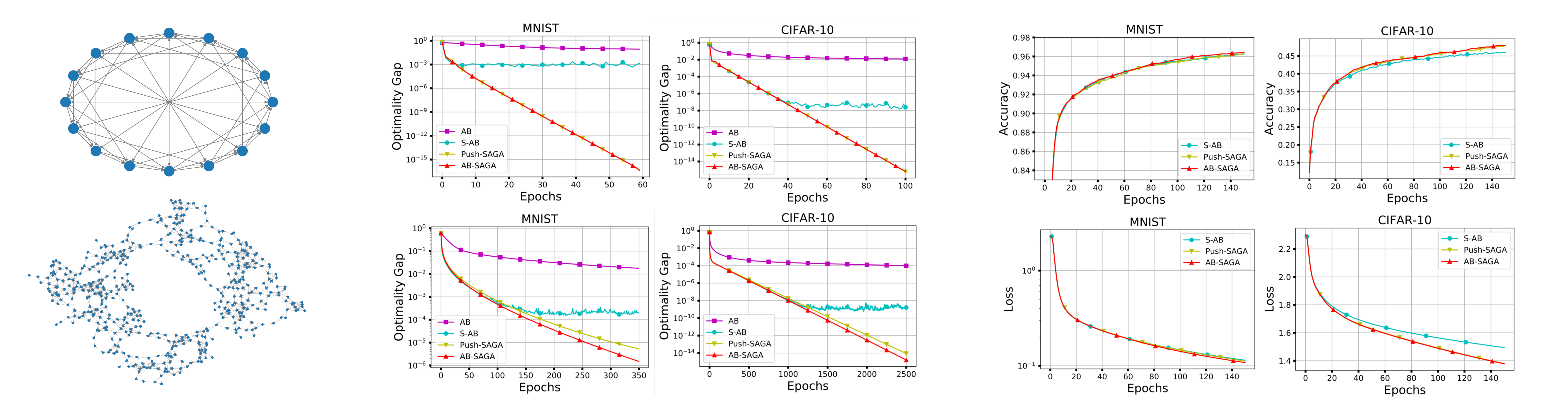}
\caption{(Left) Directed exponential graph (top) with~${n=16}$ nodes and directed geometric graph (bottom) with~${n=500}$ nodes. (Center) Optimality gap for logistic regression classifier trained over directed exponential graph (top) and directed geometric graph (bottom); (right) Test accuracy and training loss for neural networks trained over a geometric graph.}
\label{plt}
\end{figure*}

We note that~\eqref{sec},~\eqref{trd}~and~\eqref{frth}~are valid for a range of step-size, and communication rounds~$c$~and~$d$, when their right hand sides are positive. To this aim, we first fix the elements of~$\bds{\delta}$ independent of the step-size and then find the bounds on~$\alpha$. We set~${\delta_1=1}$,~${\delta_2=\frac{64 \tau_2 \kappa^2}{\ul{\pi}_r}}$,~${\delta_3 = \frac{130 \tau_2 \kappa^2 M}{m \ul{\pi}_r} }$ and~${\delta_4 = \frac{40000 n \kappa^2 M \ul{\pi}_r^{-1} \ul{\pi}_c^{-1}}{m \tau_1 (1-\sigma_B^{2d})}}$, when~${\sigma_B^{d}<\frac{\bpi_r^\top \bpi_c}{201 \kappa} \sqrt{\frac{m}{n M h_c}}}$ and ${\tau_1 := 1 - \tfrac{40000 n \sigma_B^{2d} \kappa^2 M h_c}{m(1-\sigma_B^{2d})(\bpi_r^\top \bpi_c)^2}, \tau_2 := 1 + \tfrac{40000 n \kappa^2 M h_c }{(\bpi_r^\top \bpi_c)^2 \tau_1 m (1-\sigma_B^{2d})}}$. It can be verified that the for these values of~$\bds{\delta}$, the right hand sides of~\eqref{sec},~\eqref{trd}~and~\eqref{frth} are positive. We next solve for finding a range for the step-size~$\alpha$. From~\eqref{sec}, we have
\begin{align*}
     \alpha &< \frac{g_5}{2 g_6} \frac{\delta_2}{\delta_3} - \frac{g_4}{g_6} \frac{\delta_1}{\delta_3} - \frac{g_7}{g_6} \frac{\delta_4}{\delta_3}, \\
     \impliedby \alpha &\leq \frac{1}{\kappa \ell} \left( \frac{m}{135 M (\bpi_r^\top \bpi_c)}\right).
\end{align*}
Similarly, plugging these values of~$\bds{\delta}$, in~\eqref{trd} yields
\begin{align*}
    \alpha &\leq \frac{2}{M g_5} - \frac{4 \ul{\pi}_r^{-1}}{m g_5} \frac{\delta_1}{\delta_3} - \frac{4}{m g_5} \frac{\delta_2}{\delta_3},\\
    \impliedby \alpha &\leq \frac{1}{\mu} \left( \frac{1}{9 M (\bpi_r^\top \bpi_c)}\right).
\end{align*}
To find a bound on~$\alpha$~from~\eqref{frth}, we need to ensure that ${\sigma_B^{d} \leq \frac{\bpi_r^\top \bpi_c}{1265 \kappa } \sqrt{\frac{m}{n M h_c}}}$ and therefore, we have
\begin{align*}
    \alpha &\leq \frac{1}{2 g_5} - \frac{\sigma_B^{2d}}{g_5 \delta_4} \left(\frac{292 n \ul{\pi}_r^{-1} \ul{\pi}_c^{-1}}{1-\sigma_B^{2d}} \delta_1 + \frac{194 n \ul{\pi}_c^{-1}}{1-\sigma_B^{2d}} \delta_2\right) \\
    &~~~ -\frac{\sigma_B^{2d}}{g_5 \delta_4} \left( \frac{52 \ul{\pi}_c^{-1}}{1-\sigma_B^{2d}} \delta_3 \right)\\
    \impliedby \alpha &\leq \frac{1}{225\mu} \left( \frac{1}{n (\bpi_r^\top \bpi_c)}\right).
\end{align*}
Finally, we note that~\eqref{frt} has solution if we bound ${\alpha \leq \frac{1}{\mu}\cdot\frac{2}{5 n (\bpi_r^\top \bpi_c)}}$ for the first term and~${\alpha \leq \frac{1}{\kappa \ell}}$ for the rest of the terms and ensure
\begin{align*}
    \sigma_A^{2c} &< \min \left\{\frac{m (1-\sigma_A^{2c})}{51200 n M \tau_2 h_r}, \frac{m \tau_1 (1-\sigma_A^{2c}) (1-\sigma_B^{2d})}{640000 n M h_r h_c} \right\}.
\end{align*}
To simplify the bounds on~$\sigma_A^c$ and~$\sigma_B^d$, it can be verified that~${\sigma_A^{c} < \frac{m (1-\sigma_B^{2d})}{90512 n M \kappa} \sqrt{\frac{\bpi_r^\top \bpi_c}{h_r h_c}}}$ and ${\sigma_B^d < \frac{\bpi_r^\top \bpi_c}{1265 \kappa}\sqrt{\frac{m}{n M h_c}}}$ satisfies all of the above bounds. We next define a least upper bound~$\ol{\alpha}$ on the step-size,

\begin{align*}
    \ol{\alpha} &:= \min \left\{ \frac{1}{35 \ell \sqrt{h_r h_c}}, \frac{m}{288 M n \kappa \ell (\bpi_r^\top \bpi_c)}, \frac{1}{9 \mu M (\bpi_r^\top \bpi_c)} \right\}.
\end{align*}
If~$\alpha \in (0, \ol{\alpha})$, and the communication rounds,
\begin{align*}
    c &> \frac{\log \left( \frac{m (1-\sigma_B^{2d})}{90512 n M \kappa} \sqrt{\frac{\bpi_r^\top \bpi_c}{h_r h_c}} \right)}{\log \sigma_A}, \quad d > \frac{\log \left(\frac{\bpi_r^\top \bpi_c}{1265 \kappa }\sqrt{\frac{m}{n M h_c}}\right)}{\log \sigma_B};
\end{align*}
from Lemma~\ref{spec_rad_lem}, the spectral radius~$\rho(G_\alpha) \leq \gamma = 1 - \frac{\alpha \mu n \bpi_r^\top \bpi_c}{8}$. Furthermore, if~$\alpha = \ol{\alpha}$ and~$\psi:= \frac{\sqrt{h_r h_c}}{n (\bpi_r^\top \bpi_c)}$,
\begin{align*}
    \rho(G_\alpha) \leq 1 - \min \left\{ \frac{1}{35 \kappa \psi}, \frac{m}{288 \kappa^2 M}, \frac{1}{9 M} \right\}.
\end{align*}
and the theorem follows.
\qed

\section{Numerical Experiments} \label{sims}
In this section, we illustrate the performance of~$\AS$ and compare it with related methods for finite sum minimization problems distributed over directed network of nodes. 

\noindent \textbf{Logistic Regression}: We consider a binary classification problem, using logistic regression with a strongly convex regularizer, for ${N=12,\!000}$ labelled images taken from the MNIST and CIFAR-10 datasets. These images are distributed among~$n$ nodes communicating over strongly connected directed exponential and geometric graphs, see Fig.~\ref{plt} (left). We compare~$\AB,~\SAB,~\PS$ and~$\AS$ and plot their optimality gaps~$F{(\ol{\mb x}^k) - F(\mb x^*)}$ with respect to the number of epochs, where~${\ol{\mb x}^k:=\tfrac{1}{n} \sum_{i=1}^{n} \mb x^k_i}$. We note that one epoch is~$m_i$ updates for stochastic methods and a single update for~$\AB$. It can be seen, in Fig.~\ref{plt} (center), that~$\AS$ converges linearly to the optimal solution outperforming all other methods. It is significant to note that~$\PS$ converges slower because it additionally implements the iterations for the right Perron eigenvector estimation.

\noindent \textbf{Neural Networks}: Next we consider multi-class classification problem using distributed neural networks for~${N=60,\!000}$ images taken from the MNIST and CIFAR-10 datasets. Each node trains its local neural network consisting of a hidden layer with 64 neurons and a fully connected output layer with 10 neurons. We plot the training loss~$F(\ol{\mb{x}}^k)$ and test accuracy for stochastic methods:~$\SAB, \PS$ and~$\AS$ in Fig.~\ref{plt} (right) for both graphs shown in Fig.~\ref{plt} (left). It can be observed that~$\AS$ achieves a lower loss and improved test accuracy over the other methods.

\section{Conclusions} \label{conc}
This paper describes a first-order stochastic method to minimize a distributed optimization problem such that the nodes communicate over a strongly connected directed graph. We show linear convergence of proposed method~$\AS$ to the optimal solution under weaker assumptions compared to earlier work. We also quantify the directivity constant that depicts the effects of directed nature of communication network and linear speed-up of~$\AS$ as compared to its centralized counterpart. Numerical experiments illustrates the convergence guarantees for strongly convex and nonconvex neural networks.

\bibliographystyle{IEEEbib}
\bibliography{Qureshi,Khan_allPubs,sample}
\end{document}